\documentclass[%
prx,
 amsmath,amssymb,
 aps,
twocolumn,
longbibliography,
showkeys
]{revtex4-1}

\usepackage{graphicx}
\usepackage{epstopdf}
\usepackage{dcolumn}
\usepackage{bm}
\usepackage{color}
\usepackage[mathlines]{lineno}
\usepackage[breaklinks=true,colorlinks=true,linkcolor=blue,urlcolor=blue,citecolor=blue]{hyperref}
\usepackage{ragged2e}
\usepackage{enumitem}
\usepackage{bibunits}
\usepackage{amsthm}

\usepackage{diagbox}

\defaultbibliographystyle{elsarticle-num}
\defaultbibliography{biblio}

\newcommand{\lr}[1]{\left(#1\right)}

\renewcommand\thefigure{\arabic{figure}}

\begin{document}

\title[]{Optional participation only provides a narrow scope for sustaining cooperation}

\author{Khadija Khatun$^{1,2}$}
\author{Chen Shen$^3$}
\email{steven\_shen91@hotmail.com}
\author{Jun Tanimoto$^{3,1}$}
\email{tanimoto@cm.kyushu-u.ac.jp}
\affiliation{
\vspace{2mm}
\mbox{1. Interdisciplinary Graduate School of Engineering Sciences, Kyushu University, Fukuoka, 816-8580, Japan}
\mbox{2. Department of Applied Mathematics, University of Dhaka, Dhaka-1000, Bangladesh}
\mbox{3. Faculty of Engineering Sciences, Kyushu University, Kasuga-koen, Kasuga-shi, Fukuoka 816-8580, Japan}
}
\date{\today}

\begin{abstract}

Understanding how cooperation emerges in public goods games is crucial for addressing societal challenges. While optional participation can establish cooperation without identifying cooperators, it relies on specific assumptions---that individuals abstain and receive a non-negative payoff, or that non-participants cause damage to public goods---which limits our understanding of its broader role. We generalize this mechanism by considering non-participants' payoffs and their potential direct influence on public goods, allowing us to examine how various strategic motives for non-participation affect cooperation. Using replicator dynamics, we find that cooperation thrives only when non-participants are motivated by individualistic or prosocial values, with individualistic motivations yielding optimal cooperation. These findings are robust to mutation, which slightly enlarges the region where cooperation can be maintained through cyclic dominance among strategies. Our results suggest that while optional participation can benefit cooperation, its effectiveness is limited and highlights the limitations of bottom-up schemes in supporting public goods.
\end{abstract}
\maketitle

\section{Introduction}

The maintenance of public goods---such as environmental protection, national defense, and infrastructure development---relies on cooperation among individuals~\cite{kaul2003providing}. However, cooperation often breaks down due to free-riding~\cite{archetti2012game}, where individuals benefit without contributing. To mitigate free-riding and sustain cooperation, institutional mechanisms like centralized punishment~\cite{sigmund2010social,baldassarri2011centralized} and decentralized peer-enforced measures have been proposed~\cite{fehr2002altruistic,dreber2008winners,rand2009positive}. Yet, centralized systems are susceptible to corruption~\cite{lee2019social}, while decentralized methods struggle with issues like second-order free-riding---where non-punishing cooperators benefit from the punishment efforts of punishing cooperators---and antisocial punishment, which involves the unjust targeting of cooperators~\cite{milinski2008punisher,herrmann2008antisocial,rand2011evolution}. Both approaches also require identifying free riders, which is often challenging in practice. In contrast, optional participation offers a simpler and potentially more effective mechanism by introducing a dynamic balance among cooperation, defection, and non-participation~\cite{hauert2002volunteering,semmann2003volunteering,shi2020freedom}.

Optional participation, often referred to as the ``loner'' strategy~\cite{hauert2002replicator}, allows individuals to abstain from participating in public goods interactions, thereby destabilizing the mutual defection equilibrium. Because mutual cooperation offers higher benefits than opting out, this ``escape hatch'' mechanism fosters a cycle of strategic adjustments that promote cooperative behavior. Another variant assumes that non-participants not only abstain but also harm participants without personal gain~\cite{arenas2011joker,requejo2012stability}, and this spiteful behavior can also support cooperation. Given the bottom-up nature of optional participation, these scenarios represent different non-participant behaviors: individualistic behavior, where individuals focus solely on personal payoff, and spiteful behavior, where non-participants actively harm others without deriving personal benefit. Although these strategies can effectively sustain cooperation, they capture only a limited spectrum of behaviors. This narrow scope restricts our understanding of the broader role of optional participation in fostering cooperation, particularly given the diversity of motivations and behaviors observed in social contexts.

To investigate the broader role of optional participation in cooperation, we generalized non-participation in public goods games by introducing two parameters: the outside payoff ($\alpha$), representing the net effect of abstaining and can be adjusted as an incentive or penalty, and the direct influence ($\beta$), reflecting the impact of non-participants on public goods. The values and signs of these parameters can be mapped to Social Value Orientation (SVO) categories, which reflect various ~\cite{murphy2011measuring}. While SVO typically measures how individuals balance their own welfare against others’, it also captures strategic preferences, from altruistic and prosocial to individualistic and competitive~\cite{bogaert2008social}. For example, an individualistic orientation corresponds to a positive outside payoff ($\alpha > 0$) and no influence on public goods ($\beta = 0$), reflecting a focus on personal payoff, while a sadistic orientation corresponds to a negative influence on public goods ($\beta < 0$) and a neutral outside payoff ($\alpha = 0$), indicating a tendency to harm others without personal gain. Alongside these non-participants characterized by various SVO categories, our model also includes cooperators and defectors. Cooperators contribute to the public good, while defectors do not. The total contributions are enhanced and shared equally among cooperators and defectors, with the influence ($\beta$) of non-participants also equally distributed.

Using replicator dynamics, we find that optional participation creates a limited pathway for supporting cooperation. Specifically, only non-participants motivated by individualistic or prosocial orientations can help sustain cooperation by allowing for a cyclic dominance pattern among strategies. In this pattern, cooperation gives way to defection, defection to non-participation, and non-participation back to cooperation, thus enabling coexistence among strategies. The cyclic dominance, characterized by heteroclinic cycling, is observed exclusively in the traditional loners model. In contrast, this cyclic dynamic is replaced by a stable fixed point in our generalized model, resulting in a steady-state distribution among strategies rather than ongoing cycles. The inclusion of a high mutation rate slightly expands the region where cyclic dominance supports cooperation, leading to either stable coexistence or endless cycles among strategies, depending on the sign of the direct influence on public goods. These results extend our understanding of optional participation by highlighting its limited effectiveness in promoting cooperation. They suggest that the impact of optional participation depends significantly on participants' strategic motives, and they emphasize the potential constraints of bottom-up mechanisms in sustaining public goods.

\section{Model} In an infinitely large, well-mixed population we consider a sample group of $M>2$ players randomly chosen in a public good game (PGG). In the PGG,  cooperators ($C$) contribute a fixed cost $c$ to a common pool, while defectors ($D$) contribute nothing. The total contribution is multiplied by an enhancement factor $r$ ($1<r<M$) and then equally divided among all participants, regardless of their initial contribution. While the group benefits most if all individuals contribute to the common pool, defection is the best individual strategy to maximize personal gains in the context. With the lure of defection, cooperation tends to break down in public goods game.

\begin{figure}[!t]
\includegraphics[width=1.05 \linewidth]{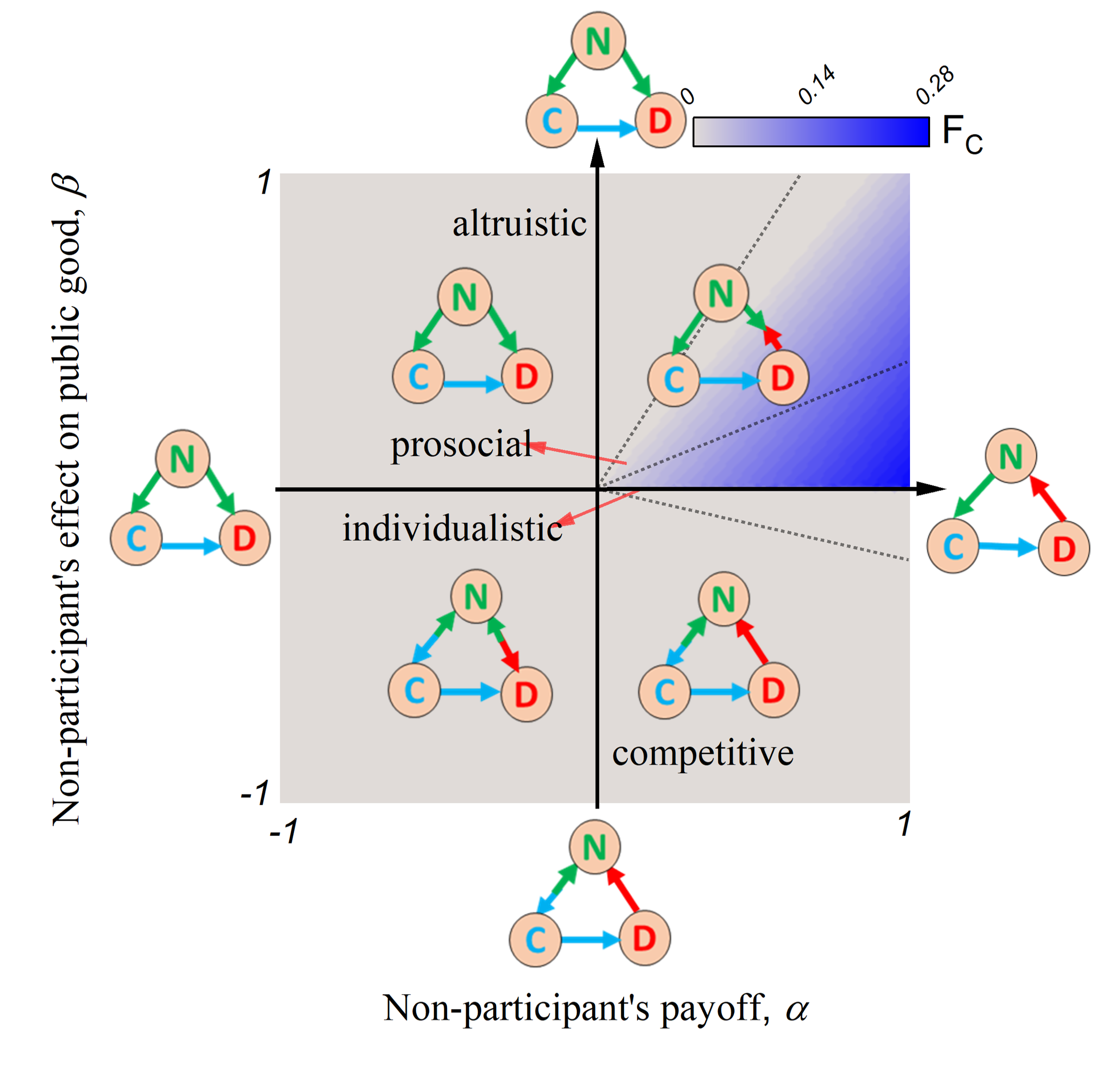}
    \caption{Non-participation narrowly sustains cooperation. Cooperation only survives if $\alpha, \beta>0$, when the non-participant's motivation is individualistic and prosocial.  Optimal cooperation in the individualistic behavior; While an increase in $\alpha$ promotes cooperation, an increase in $\beta$ inhibits it. The figure depicts a phase diagram of cooperation fraction ($F_C$) and invasion dynamics among cooperation ($C$), defection ($D$), and non-participation ($N$) strategies, as influenced by non-participants' payoff ($\alpha$) and their impact on the public good ($\beta$). Dotted lines separate the individual strategic motivation according to SVO.  The parameter values are fixed at $M=5$, $r=3$, and chose a real mutation $\mu=10^{-8}$ to avoid computational error. All subsequent plots will use the same $M$ and $r$.  }   
    \label{fig1}
\end{figure}

Building on the fundamental optional PGG~\cite{hauert2002replicator}, in which participants can opt out of PGG with a constant payoff but do nothing to the public good, we have generalized the behavior of non-participates ($N$). Particularly, non-participants, who can both receive a payoff $\alpha$ ($-1 \leq \alpha \leq 1$) and simultaneously provide a payoff  $\beta$ ($-1 \leq \beta \leq 1$), which is received by both cooperators and defectors equally. The behavioral motivations of a non-participant with the public good players can be defined by $\theta = \arctan\lr{\frac{\beta}{\alpha}}$ using SVO, where $\theta$ represents the angle of orientation~\cite{murphy2011measuring}. Therefore, the behavioral motivation is 
 (1) \textit{altruistic}, if $\theta>57.15^{\circ}$, (2) \textit{prosocial}, if $22.45^{\circ}<\theta<57.15^{\circ}$, (3) \textit{individualistic}, if $-12.04^{\circ}<\theta<22.45^{\circ}$, and (4) \textit{competitive}, if $\theta<-12.04^{\circ}$(indicated in Figure~\ref{fig1}). 
 
 Players can choose to join or refuse the public good, and those who join can then cooperate or defect. Firstly, we assume that $S$ players are willing to join the public good and $M-S$ will be non-participant's. If $m$ is the number of cooperators in $S$, and the cooperation cost $c=1$, then the payoffs of a defector, cooperator, and non-participant are:
\begin{equation}
   \begin{array}{l}
        \Pi_D=\frac{rm+\beta(M-S)}{S},   \\
        \Pi_C=\Pi_D-1,\\
        \Pi_{N}=\alpha.
   \end{array}
    \label{eq01}
\end{equation}
Let  $x$, $y$, and $z$ be the fraction of cooperation, defection, and non-participation, where $x,y,z\geq 0$ and $x+y+z=1$. If the number of participants joining the public good is $S=1$, then the player have to play as a non-participant and obtain $\alpha$, which will happen with probability $z^{M-1}$.  However, the probability of  $S-1$ ($S>1$) players among $M-1$ players joining the public good will be
${M-1 \choose S-1} (1-z)^{S-1} z^{M-S}$.\\

The probability that $m$ of those players are cooperators, and $S-1-m$ defectors is
${S-1 \choose m}\lr{\frac{x}{x+y}}^{m} \lr{\frac{y}{x+y}}^{S-1-m}$.\\ 

Since the payoff of a defector is $\frac{rm}{S}$, so the average payoff of a defector in a group of $S$ players($S=2,.., M$) will be\\
$\frac{r}{S}\sum_{m=0}^{S-1} {S-1\choose m}\lr{\frac{x}{x+y}}^{m} \lr{\frac{y}{x+y}}^{S-1-m} m\\=\frac{r}{S} (S-1)\frac{x}{1-z}$.\\

Therefore, the average payoff of a defector in the group of $M$ players is:\\ 
$\Pi_D=\alpha  z^{M-1}+\sum_{S=1}^{M} {M-1 \choose S-1} \lr{\frac{r}{S} (S-1) \frac{x}{1-z}+\beta \lr{\frac{M-S}{S}}} (1-z)^{S-1} z^{M-S}$.\\
Utilization of ${M-1 \choose S-1}={M \choose S}\frac{S}{M}$ and\\ $\sum_{S=0}^{M}S {M \choose S}(1-z)^S z^{M-S}=M (1-z)$, yields
 \begin{equation}   
 \begin{aligned}
     \Pi_D= & \alpha \left. z^{M-1}+ \frac{r x}{1-z} \lr{1-\frac{1-z^M}{M(1-z)}}+\right.\\
 & \left. \beta \lr{\frac{1-z^M}{1-z}-1}. \right.
\end{aligned}
\label{eq02}
 \end{equation}
 A defector who chooses not to cooperate avoids the cooperation cost while causing a reduction of public good by a factor $\frac{r}{S}$, then the payoff change will be $\lr{1-\frac{r}{S}}$. Consequently,\\ 
 $\Pi_{D}-\Pi_{C}= \sum_{S=2}^{M}\lr{1-\frac{r}{S}} {M-1 \choose S-1}(1-z)^{S-1} z^{M-S}$.\\
 Therefore,
 \begin{equation}
 \begin{array}{l}
    \Pi_{D}-\Pi_{C}=1+(r-1) z^{M-1}- \frac{r}{M} \frac{1-z^M}{1-z}=:G(z).
 \end{array}
\label{eq03}
 \end{equation}

 So, the average payoff of a cooperator is:
 \begin{equation}
     \begin{aligned}
     \Pi_C
 =-1- &\left. (r-1) z^{M-1}+ \frac{r}{M} \frac{1-z^M}{1-z}+ \alpha z^{M-1}\right.\\
 & \left. + \frac{r x}{1-z} \lr{1-\frac{1-z^M}{M(1-z)}}+\beta \lr{\frac{1-z^M}{1-z}-1} \right..
        \end{aligned}
        \label{eq04}
 \end{equation}
The average payoff in the whole population is:
 \begin{equation}
     \begin{array}{l}
     \overline{\Pi}=x*\Pi_{C}+y*\Pi_{D}+z*\Pi_{N}.
        
     \end{array}
     \label{eq05}
 \end{equation}

Therefore, the replicator-mutator dynamics with mutation $\mu$ are:
\begin{equation}
\begin{array}{l}
\dot{x}=x \lr{\Pi_{C}-\overline{\Pi}}-\mu x+\frac{\mu (1-x)}{2},\\
\dot{y}=y \lr{\Pi_{D}-\overline{\Pi}}-\mu y+\frac{\mu (1-y)}{2},\\
\dot{z}=z \lr{\Pi_{N}-\overline{\Pi}}-\mu z+\frac{\mu (1-z)}{2}.
\end{array}
\label{eq06}
\end{equation}

Detailed explanations of the equilibria and their stability of the replicator dynamics have been given in the Appendix. 

\begin{figure*}[!t]
\includegraphics[width=1\linewidth]{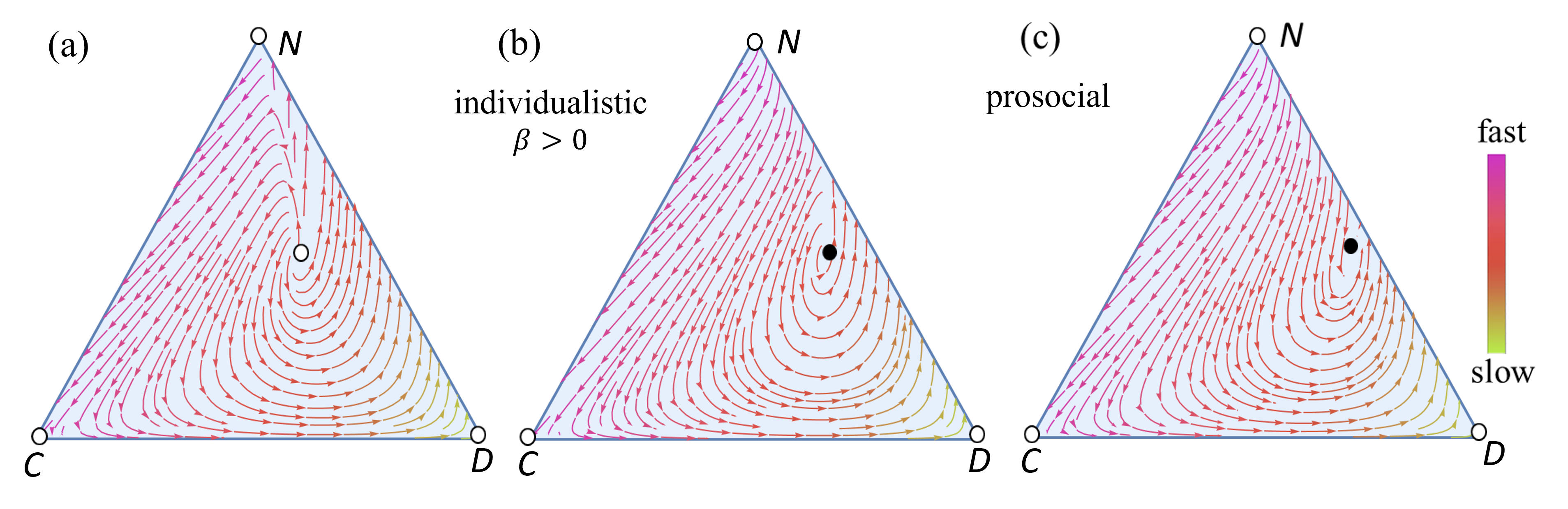}
    \caption{Cooperation survives. (a) Hetero-clinic cycle among cooperation, defection, and non-participant when $\alpha>0$ and $\beta=0$, and stable coexistence of those three when $\alpha, \beta>0$ (b)individualistic and (c) prosocial. With the increase of benefits to the public good equilibrium goes near to the coexistence of defection and non-participation.   The parameter values are fixed at (a) $\alpha=0.5$, $\beta=0$, (b) $\alpha=0.5$, $\beta=0.2$, and (c) $\alpha=0.8$, $\beta=0.7$, and $\mu=10^{-8}$ for all. Solid black dots are stable and whites are unstable equilibrium points.
    }   
    \label{fig2}
\end{figure*}

\begin{figure*}[!t]
\includegraphics[width=1.0\linewidth]{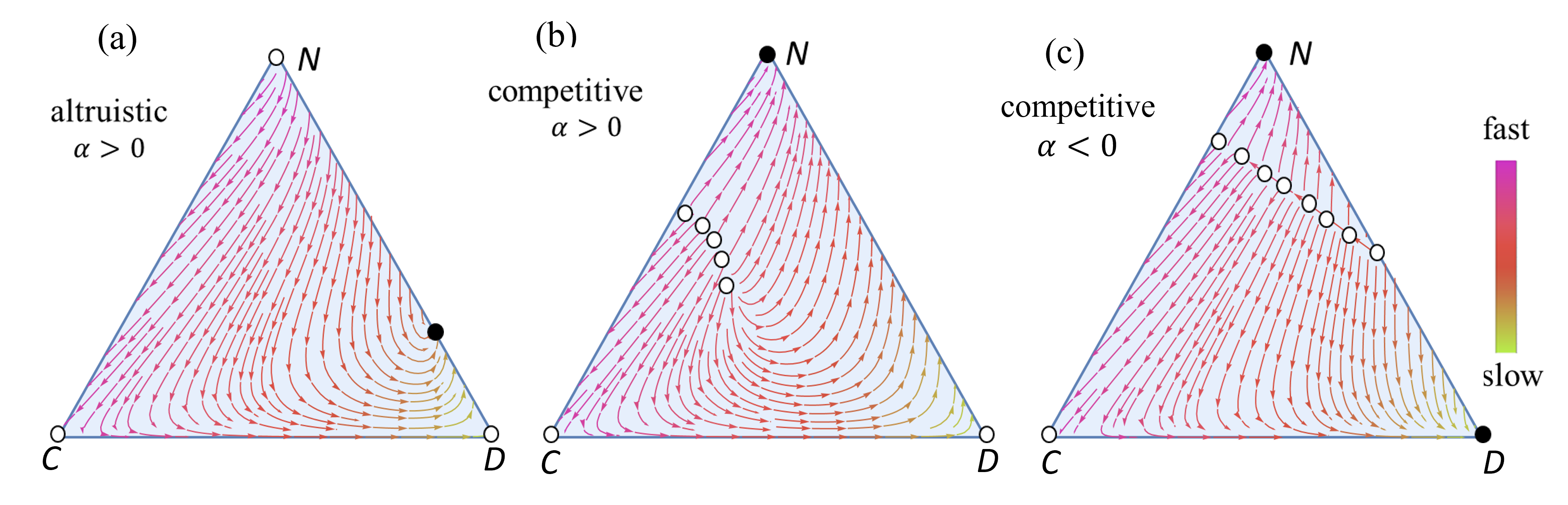}
    \caption{Cooperation can not survive. (a) Cooperation overtook non-participation but destabilization of defection is not non-participation dominant in the coexistence of defection and non-participation when $\alpha<\beta$ (in altruism). (b) Cooperation can not invade non-participation fully; mono-stable non-participation when $\alpha>0$ and $\beta<0$, and (c) neither cooperation invade non-participation fully nor destabilization of defection is non-participation dominant; bi-stable between defection and non-participation when $\alpha, \beta<0$ (some potion of individualistic and competitive motivations).  The parameter values are fixed at $\mu=10^{-8}$, (a) $\alpha=0.4$, $\beta=0.9$, (b) $\alpha=0.5$, $\beta=-0.8$, and (c) $\alpha=-0.5$, $\beta=-0.8$. Solid black dots are stable and whites are unstable equilibrium points.}
   
    \label{fig3}
\end{figure*}
\begin{figure}[!t]
\includegraphics[width=1.0\linewidth]{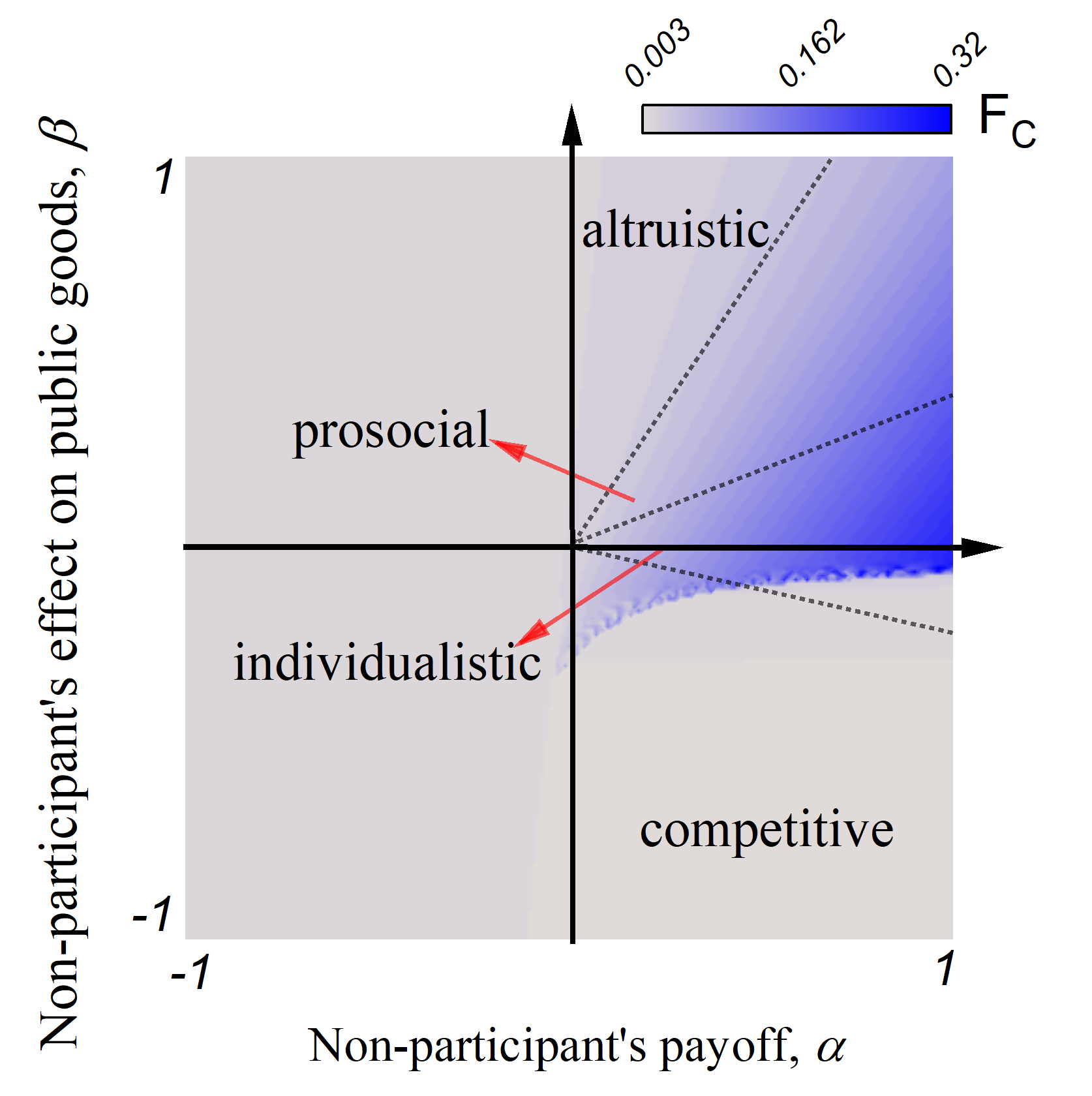}
    \caption{Large mutation does not alter the results but slightly enlarges cooperation region. Only for non-negative $\alpha$ mutation can maintain cooperation.  Mutation rate is, $\mu=10^{-2}$.   }   
    \label{fig4}
\end{figure}
\begin{figure*}[!t]
\includegraphics[width=1.0\linewidth]{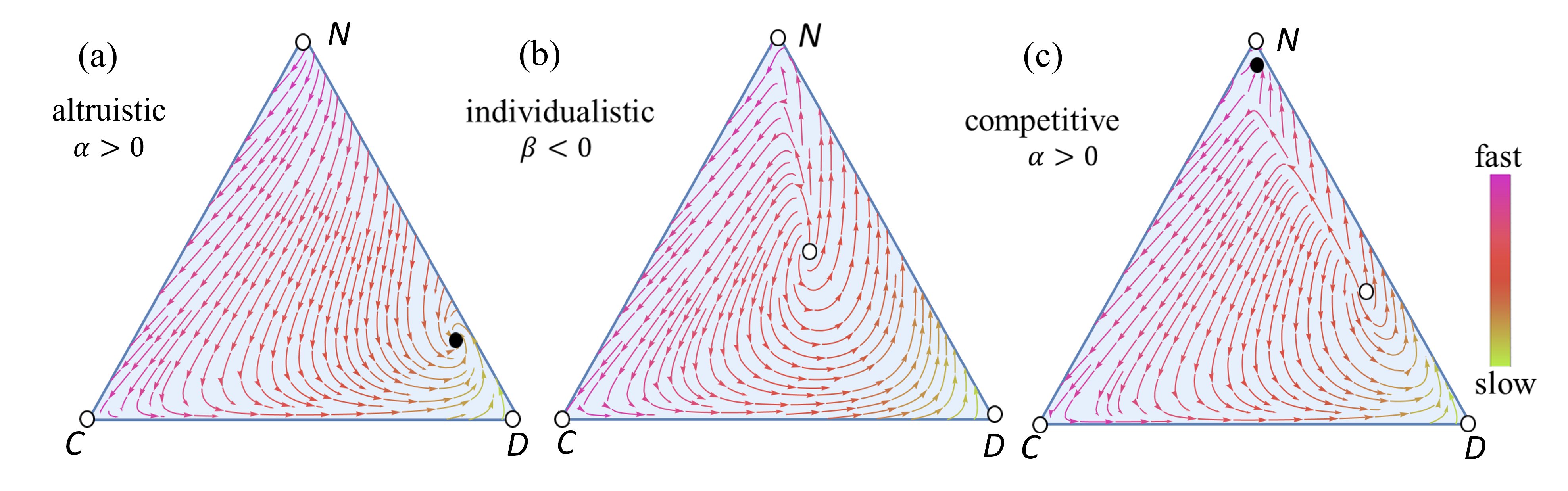}
    \caption{ (a) A stable coexistence of cooperation, defection, and non-participating in the altruistic,  (b) hetero-clinic cycle in the individualistic motivation when non-participants incur damage to the public good, and (c) stable coexistence of three strategies after oscillation in competitive motivation. The parameter values are fixed at (a) $\alpha=0.3$, $\beta=0.8$, (b) $\alpha=0.8$, $\beta=-0.05$, and (c) $\alpha=0.125$, $\beta=-0.25$, and $\mu=10^{-2}$ for all. Solid black dots are stable and whites are unstable equilibrium points. }   
    \label{fig5}
\end{figure*}

\section{Results}
Figure~\ref{fig1} shows a phase diagram of the fraction of cooperators ($F_C$) as a function of the non-participants' payoff $\alpha$ and their impact on public goods $\beta$. Specifically, when $\alpha > 0$ and $\beta = 0$, our model reduces to the loner model~\cite{hauert2002replicator}. Conversely, when $\alpha = 0$ and $\beta < 0$, it corresponds to a variant of the jokers model~\cite{arenas2011joker}, with a slight modification: the minimum required number of public goods participants is greater than one. As expected, we confirmed the positive effect of loners on cooperation, as extensively documented in previous studies~\cite{hauert2002replicator, mathew2009does, rand2011evolution}. However, under rare mutation scenarios, we did not observe cooperation in the joker model, because of model differences. 
 
Beyond these two specific scenarios, our analysis shows that non-participation's role in sustaining cooperation is quite narrow, as illustrated in Figure~\ref{fig1}. The results demonstrate that cooperation can only persist when  $\alpha>0$ and $\beta \geq 0$ correspond to non-participants' individualistic or prosocial motivations. Notably, cooperation is optimum under individualistic behavior, specifically when non-participants receive substantial self-payoff ($\alpha \to 1$) while contributing minimally to the public good ($\beta \to 0$).
 
To further understand these phenomena, we examined the invasion dynamics between each pair of strategies, as illustrated in the triangular diagrams in Figure~\ref{fig1} (detailed theoretical analysis is provided in the Appendix; see Figures~\ref{fig6} and \ref{fig7}). We found that the relationship between cooperation and defection remains consistent regardless of the non-participants' payoff and their impact on public goods: defection always invades cooperation. However, the dynamics involving non-participation and the other two strategies are more complex and deviate from the straightforward hierarchical dominance patterns typically observed in loner or joker models. Specifically, when $\alpha > 0$ and $\beta \geq 0$, cooperation invades non-participation while coexisting with defection; when $\alpha \geq 0$ and $\beta < 0$, non-participation invades defection and is bistable with cooperation; when $\alpha < 0$ and $\beta < 0$, non-participation is bistable with both cooperation and defection; and when $\alpha \leq 0$ and $\beta \geq 0$, non-participation is invaded by both cooperation and defection.  
 
In our generalized model, defectors consistently invade cooperators regardless of the non-participants' payoff ($\alpha$) or their impact on public goods ($\beta$). This means that cooperation can only survive through a specific pathway involving non-participation: defectors must be destabilized by non-participants, who, in turn, must be destabilized by cooperators. This sequence is feasible in most regions except where the outside payoff $\alpha$ is negative and the influence on public goods $\beta$ is positive. However, cooperation is only observed under conditions where $\alpha$ is positive and $\beta$ is non-negative, as shown in Figure~\ref{fig1}. In this case, cooperators can fully dominate non-participants (see Figure~\ref{fig6}(a)), and non-participants can coexist with defectors (see Figure~\ref{fig7}(a)). For cooperation to survive in the whole system, it is necessary that non-participants overtake defectors when coexisting with them, which generally occurs when $\alpha$ is greater than $\beta$. This relaxes the strict condition in the traditional loners model, where non-participation must fully dominate defectors. This strict condition leads to a heteroclinic cycle of cooperation, defection, and non-participation in the loners model, as confirmed in Figure~\ref{fig2}(a). In contrast, allowing positive $\beta$ in our model leads to the stable coexistence of these strategies (Figures~\ref{fig2}(b) and \ref{fig2}(c)). Conversely, when the direct negative influence of non-participation on public goods ($\beta$) outweighs the non-participants' outside payoff ($\alpha$), cooperation disappears, and non-participation coexists with defection (Figure~\ref{fig3}(a)).

Although the conditions where $\alpha \geq 0$ and $\beta < 0$, and where $\alpha < 0$ and $\beta < 0$ might theoretically support the survival of cooperation, we do not observe cooperation in these regions. In the former case, trajectories either lead directly to non-participation or involve a sequence where defection invades cooperation, followed by a shift to non-participation (see Figure~\ref{fig3}(b)). In the latter case, trajectories proceed directly to non-participation or defection after initially invading cooperation (see Figure~\ref{fig3}(c)). Previous studies have demonstrated that moderate mutation rates in repeated games can enhance cooperation\cite{tkadlec2023mutation}. Building on this insight, we investigated whether larger mutation rates might further affect cooperative dynamics. Specifically, we examined how cooperation might persist in these parameter regions under conditions of high mutation.
 
A larger mutation rate ($\mu=0.01$) slightly expands the region in which cooperation can survive, as shown in Figure~\ref{fig4}. This expansion does not include regions with negative $\alpha$, where the emergence of cooperation would depend solely on mutation. Instead, the expansion appears in two specific parameter regions: when $0 < \alpha < \beta$ and when $\alpha \geq 0$ and $\beta < 0$. 
In the first parameter region, characterized by altruistic motivation, the system's behavior partially mirrors that observed under rare mutation. While the pathways from non-participation to cooperation remain unchanged from the rare mutation scenario (Figure \ref{fig3}(a)), the introduction of larger mutation rates induces a subtle shift in the defection-nonparticipation equilibrium, which enables a stable three-strategy coexistence state involving cooperation, defection, and non-participation (Figure \ref{fig5}(a)). In the second parameter region, where non-participants inflict minimal damage to the public good ($\beta$ approaching $0$ from below), we observe distinct dynamics in both individualistic and competitive motivational contexts. Here, cooperation can persist through two mechanisms: either via heteroclinic cycles among non-participation, cooperation, and defection (Figure \ref{fig5}(b)) or through a stable three-strategy coexistence (Figure \ref{fig5}(c)). Notably, despite requiring more than one public good participant, our findings at $\alpha=0$ and $\beta<0$ qualitatively replicate the stable coexistence observed in the Joker model under high mutation rates.

Although large mutation rates allow us to observe sustained cooperation in regions where non-participants have positive payoffs and negatively influence public goods, we fail to observe cooperation in regions with negative $\alpha$ and $\beta$, even when the mutation rate is large. Despite invasion direction analysis for strategy pairs indicating the possibility of cooperation, survival in this region may require additional conditions---such as population structure or heterogeneous mutation rates among strategies---which warrants further investigation.

\section{Discussion}

Our study aimed to explore the broader role of optional participation in fostering cooperation in one-shot, anonymous public goods games. To this end, we generalized the behavior of non-participants by introducing two simple parameters: the outside payoff of non-participation and their direct influence on public goods. We found that optional participation sustains cooperation only when non-participants are individualistically or prosocially motivated, with individualistic non-participation leading to optimal cooperation. Considering mutation did not significantly alter these results, except by slightly expanding the region in which cooperation is supported through cyclic dominance dynamics. These findings suggest that while optional participation can benefit cooperation, its effectiveness is limited.

The optional participation mechanism, primarily represented by the loner model~\cite{hauert2002replicator,hauert2002volunteering} and its variant, the jokers model~\cite{arenas2011joker,requejo2012stability}, is an oversimplified approach that limits its applicability. In these models, individuals can choose to opt out of the public goods game. In the loner model, non-participants (loners) receive a fixed payoff without directly affecting the public good, while in the jokers model, non-participants receive no payoff but actively harm others. These simplified assumptions produce cycles among cooperation, defection, and non-participation, which can support cooperative behavior, though this effect in the jokers model only emerges under high mutation rates. While the loner and jokers models capture specific instances of optional participation, they restrict the possible range of non-participant behaviors, especially given the bottom-up nature of optional participation schemes. For example, consider a bridge connecting two villages across a river: individuals may choose to use the bridge or to opt out by using their own boats. The impact of opting out varies depending on individual actions; spiteful individuals might damage the bridge despite not using it, while prosocial individuals might contribute to its maintenance even if they do not use it. Additionally, the jokers model assumes that opting out incurs no cost, which is often unrealistic; in practice, non-participation frequently involves costs, so the payoff for non-participants may not always be positive or zero, as assumed in these models.

By generalizing the optional participation model to incorporate both the outside payoff of non-participation and its direct impact on public goods, we demonstrated that cooperation can persist under less restrictive conditions than those required in the loner and jokers models. These earlier models necessitate strict cyclic dominance among strategies, where each strategy must strictly dominate another in a closed loop, resulting in a continuous cycle among cooperation, defection, and non-participation. In our generalized model, however, cooperation can survive if two key conditions are met. First, non-participation must destabilize defection, preventing defectors from fully dominating the population. Second, under these circumstances, cooperation must either completely eliminate non-participation or establish itself as the more stable equilibrium in bi-stable states where cooperation and non-participation compete. In the former case, cooperation stably coexists with both defection and non-participation, achieving a state that is resilient to mutation. In the latter case, cycles among cooperation, defection, and non-participation emerge, but only under conditions of high mutation rates.

Our results have significant socioeconomic implications, suggesting that bottom-up schemes for public goods, which empower individual decision-making, may not always foster cooperation. Policymakers must consider the direct, active influence of non-participants on public goods. For example, in environmental conservation, individuals who opt out may harm the environment through littering, illegal dumping, or poaching. In public health, those who refuse vaccinations might spread misinformation or ignore quarantine, increasing disease transmission. In communal infrastructure, farmers who don’t contribute may divert water illegally or sabotage shared systems. Sustaining cooperation requires balancing power relations by adjusting the outside payoff of non-participation and countering its influence on public goods. This could include penalties for harmful actions or incentives for behavior that protects public goods, ensuring the system functions even when individuals can opt-out.

In sum, by generalizing the actions of non-participation, we find that optional participation has a limited role in the survival of cooperation. Despite this, future work should consider the heterogeneous influences of non-participation on public goods and explore more realistic scenarios like repeated games~\cite{schmid2021unified}, the costly punishment problem~\cite{rand2011evolution}, or social networks~\cite{civilini2024explosive}. Such studies could reveal that aspects of human cooperativeness, obscured by simplified theoretical assumptions, emerge only in complex real-world contexts, whose forces have shaped our decision-making abilities.
 
\paragraph*{Acknowledgements.} 
We also acknowledge support from (i) Japanese Government (MEXT) scholarship, Japan, ( grant No. 222143) awarded to K.\,K., (ii) JSPS KAKENHI (grant no. JP 23H03499) to C.\,S., and (iii) the grant-in-Aid for Scientific Research from JSPS, Japan, KAKENHI (grant No. JP 23H03499) awarded to J.\,T. 
\paragraph*{Author contributions.} 
C.\,S. and J.\,T. conceived research. K.\,K. and C.\, S. performed analytical analysis. All co-authors discussed the results and wrote the manuscript.
\paragraph*{Data statements and datasets in repositories} This is a theoretical study, no data was used in the analysis.

\paragraph*{Conflict of interest.} Authors declare no conflict of interest.


\section{Appendix}
\renewcommand\thefigure{A\arabic{figure}}
\setcounter{figure}{0} 
\renewcommand\theequation{A\arabic{equation}}
\setcounter{equation}{0}  
\subsection{Stationary Points and their Stability Analysis}
The replicator dynamic Eq.~\ref{eq06} without mutation($\mu=0$) has at most six stationary points: $E_1=\lr{1,0,0}$, $E_2=\lr{0,1,0}$, $E_3=\lr{0,0,1}$, $E_4=\lr{0,\frac{\beta}{\alpha+\beta},\frac{\alpha}{\alpha+\beta}}_{\alpha, \beta > 0 or  \alpha, \beta<0}$, $E_5=\lr{\frac{\beta}{1-r+\alpha+\beta},0,\frac{1-r+\alpha}{1-r+\alpha+\beta}}_{-1\leq \alpha \leq 1, -1 \leq \beta <0, and r>2}$ and $E_6=\lr{x^{*},1-x^{*}-z^{*},z^{*}}$ (here $x^{*}, z^{*}$, since the exact form is quite complicated). Due to the complexity of the analysis, we use a side-wise invasion relation to analyze the stability of the equilibrium points.\\

On $C-D$ side $x+y=1$ and $z=0$, 
\begin{equation}
\begin{array}{l}
     \dot{x}=x\lr{1-x}\lr{\Pi_{C}-\Pi_{D}}=x \lr{1-x}\lr{\frac{r}{M}-1}.
\end{array}
\label{eq07}
\end{equation}
There are only boundary equilibria i.e. defection and cooperation. Since $r<M$, the rate will be negative so that defection will be the final stage in the evolution of cooperation and defection.\\

\begin{figure}[!t]
\includegraphics[width=1\linewidth]{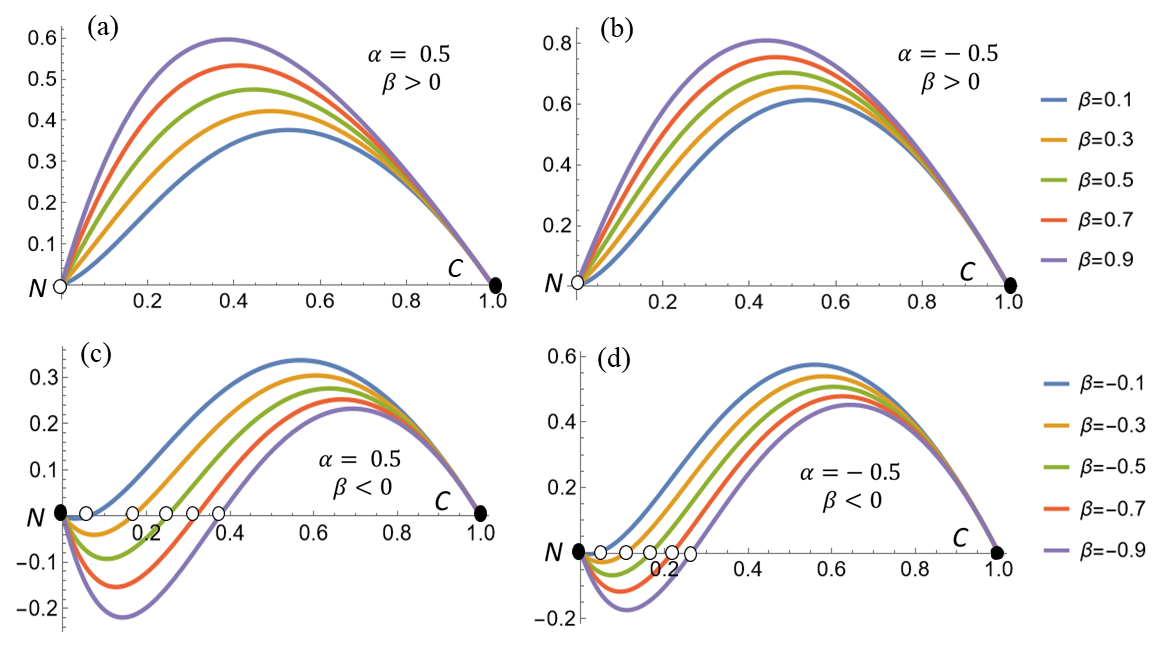}
    \caption{Stability of the equilibria on the $C-N$ side. When $\beta>0$ and (a) $\alpha>0$  or (b) $\alpha<0$ : Cooperation is stable. When $\beta<0$ and (c) $\alpha>0$  or (d) $\alpha<0$ : Bi-stability between cooperation and non-participant; most possible outcome is cooperation if $\beta$ near to $0$. The black dot is for stable points and the white dot is for unstable ones.}   
    \label{fig6}
\end{figure}

On $C-N$ side, where $x+z=1$ and $y=0$, 
\begin{equation}
\begin{array}{l}
     \dot{x}=x\lr{1-x}\lr{\Pi_{C}-\Pi_{N}}\\
     =-\lr{-1 + (1 - x)^M + x} \lr{x \lr{-1 + r - \alpha- \beta} + \beta}.
\end{array}
\label{eq08}
\end{equation}
On this side, an intermediate equilibrium ($C+N$), $\lr{\frac{\beta}{1-r+\alpha+\beta},\frac{1-r+\alpha}{1-r+\alpha+\beta}}$ is possible, when $-1\leq \alpha \leq 1$, $-1 \leq \beta <0$, and $r>2$  with boundary equilibria- cooperation and non-participation.  Cooperation will be stable for $\beta>0$ and either $\alpha>0$ or $\alpha<0$ (Figures \ref{fig6}(a) and \ref{fig6}(b)). Bi-stable between cooperation and non-participant when $\beta <0$ and either $\alpha>0$ or $\alpha<0$ (Figures \ref{fig6}(c) and \ref{fig6}(d)). Cooperation is the most possible stable outcome if $\beta \to 0$. \\

\begin{figure}[!t]
\includegraphics[width=1\linewidth]{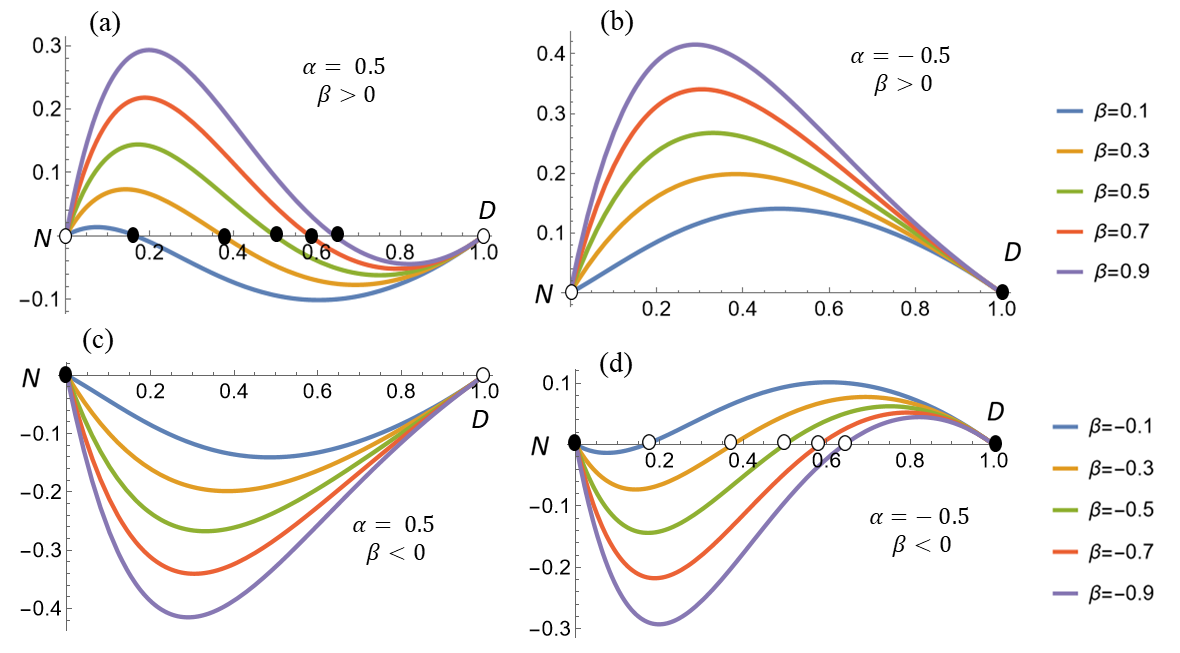}
    \caption{Stability of the equilibria on the $D-N$ side. (a) if $\alpha,\beta>0$:  Mixed equilibrium ($D+N$) will be stable; non-participation is dominant over defection if $\alpha>\beta$ . If (b) $\alpha \leq  0$, and $\beta \geq 0$: Defection is stable. (c) when $\alpha > 0$, $\beta < 0$: non-participant is mono-stable and (d) when $\alpha < 0$, $\beta < 0$: Bi-stability between defection and non-participant.  The black dot is for stable points and the white dot is for unstable ones.}   
    \label{fig7}
\end{figure}

On $D-N$ side $y+z=1$ and $x=0$,
\begin{equation}
\begin{array}{l}
     \dot{y}=y\lr{1-y}\lr{\Pi_{D}-\Pi_{N}}\\=
     (-1 + (1 - y)^M + y) (-\beta + y (\alpha + \beta)).
\end{array}
\label{eq09}
\end{equation}
There is an intermediate equilibrium ($D+N$), $\lr{\frac{\beta}{\alpha+\beta},\frac{\alpha}{\alpha+\beta}}$, for $\alpha, \beta > 0$, and $\alpha, \beta < 0$. Looking at Figure~\ref{fig7}(a), for $\alpha,\beta>0$,  the mixed equilibrium will be stable. This mixed equilibrium varies with the relative relation of $\alpha$ and $\beta$, if $\alpha>\beta$ then non-participation is dominant over defection otherwise defection is dominant. Defection is stable when $\alpha \leq 0$ and $\beta\geq 0$ (see Figure~\ref{fig7} (b)),  non-participation is stale when $\alpha>0$ and $\beta <0$ (Figure~\ref{fig7} (c)) and bi-stable between defection and non-participant, if $\alpha, \beta \leq 0$ (Figure~\ref{fig7}(d)). 

Therefore, in $C-D-N$ simplex,
when $\alpha, \beta>0$, the non-participant is invaded by cooperation on the $C-N$ side (Figure \ref{fig6}(a)). In that case, defection converges to the mix of defection and non-participant on the $D-N$ side, and the mixed equilibrium ($D+N$) varies with $\alpha$ and $\beta$ (Figure \ref{fig7}(a)). In the presence of non-participants, cooperation is preferable if non-participants exceed a threshold of $0.5$.  When $\alpha>\beta>0$, the proportion of non-participants exceeds a threshold of $0.5$, fostering cooperation, whereas when $0<\alpha <\beta$, the proportion of non-participants falls below $0.5$, insufficient to promote cooperation.   Hence, the mixed equilibrium of cooperation, defection, and non-participant is stable if $\alpha>\beta>0$, and the mixed equilibrium of defection and non-participant is stable if $0 < \alpha < \beta$.

When $\alpha>0$ and $\beta=0$, the stability analysis is inconclusive in this method, to analyze the dynamics we transfer the replicator dynamic into the Hamiltonian system. Let, $f=\frac{x}{x+y}$, the fraction of cooperation among the individuals participating in the public good game, then we get, $\dot{f}=\frac{y\dot{x}-x \dot{y}}{(x+y)^2}=\frac{xy}{(x+y)^2} \lr{\Pi_{C}-\Pi_{D}}$, and $0\leq f\leq 1$. Therefore, the Eq.\ref{eq06} (at $\mu=0$) will be,  
 \begin{equation}
\begin{array}{l}
     \dot{f}=f (1-f)\lr{\Pi_{C}-\Pi_{D}}=- f (1-f) G(z),\\
     \dot{z}=z (1-z)\lr{1-z^{M-1}} \lr{\alpha- (r-1) f-\frac{\beta}{1-z}}.
\end{array}
\label{eq10}
\end{equation}
The division of the right-hand side by the fraction $f(1-f) z (1-z) (1-z^{M-1})$, does not change the direction of the orbits on the unit square $(0,1)^2$. At $\beta=0$, which yields,
 \begin{equation}
\begin{array}{l}
     \dot{f}= -\frac{G(z)}{z (1-z) (1-z^{M-1})}=u(z),\\
     
     \dot{z}= \frac{\alpha- (r-1) f}{f (1-f)}=w(f).
\end{array}
\label{eq11}
\end{equation}
Let, $H=U+W$, where $U(z)$, and $W(f)$ are the anti-derivatives of $u(z)$, and $w(f)$: 
\begin{equation}
\begin{array}{l}
     U(z)=(1-\frac{r}{M}) \log{z}+(\frac{r}{2}-1) \log{(1-z)}+R(z).
\end{array}
\label{eq12}
\end{equation}
\begin{equation}
\begin{array}{l}
     W(f)=\alpha  \log{f}+(r-1-\alpha) \log{(1-f)}.
\end{array}
\label{eq13}
\end{equation}
$R(z)$ bounded on $[0,1]$, we get the Hamiltonian system 
\begin{equation}
\begin{array}{l}
\dot{f}=-\frac{\partial H}{\partial z}, \\
\dot{z}=\frac{\partial H}{\partial f}.
\end{array}
\label{eq14}
\end{equation}
Near the boundary of $[0,1]^2$, all level sets of $H$ are closed, since $H\to -\infty$ uniformly, when $U(z) \to -\infty$ for $z\to 0, 1$ if $2<r< M$ and $W(f) \to -\infty$ for $f\to 0, 1$ if $\alpha<r-1$. Since the Hamiltonian system is conservative, all level sets within $[0,1]^2$ are closed. So, at  $\alpha>0$ and $\beta=0$, there is a hetero-clinic Rock-Scissors-Paper type of cycling among cooperation, defection, and non-participation.  

If $\alpha>0$ and $\beta<0$, bi-stable is shown between non-participant and cooperation on the $C-N$ side, where cooperation is invaded by defection on the $C-D$ side and non-participant is stable on the $D-N$ side. Hence, $E_3$ is a stable point.

When $\alpha, \beta \leq 0$, and $\alpha= \beta \neq  0$, bi-stability is shown between non-participant and cooperation on the $C-N$ side and non-participant and defection on the $D-N$ side, cooperation can be invaded by the defection on the $C-D$ side. Hence, $E_2$ or $E_3$ is the final bi-stable state.

Finally, when $\alpha <0$, and $\beta >0$, the non-participant is invaded by both cooperation on the $C-N$ side and defection on the $D-N$ side, and the cooperation is invaded by the defection on the $C-D$ side. Hence, $E_2$ is a stable point. 
\bibliography{Bibliography}
\bibliographystyle{elsarticle-num}

\end{document}